\documentclass[a4paper,oneside]{article}
\usepackage[top=4cm,bottom=3.5cm,left=3.5cm,right=3.5cm]{geometry}
\usepackage[T1]{fontenc}
\usepackage[utf8]{inputenc}
\usepackage[UKenglish]{babel}
\usepackage{amsthm}
\usepackage{amsmath,amssymb}
\usepackage{nicefrac}
\usepackage{mathtools}
\usepackage{mathdots}
\usepackage{hyperref}
\usepackage{centernot}
\usepackage[inline]{enumitem}
\newlist{choices}{enumerate*}{1}
\setlist[choices]{itemjoin = \hspace{1cm}, label=\arabic*.}

\usepackage{thmtools}
\usepackage{cleveref}
\usepackage{nameref}
\Crefname{rem}{Remark}{Remarks}
\Crefname{alphthm}{Theorem}{Theorems}
\Crefname{thm}{Theorem}{Theorems}
\Crefname{equation}{Equation}{Equations}

\newcommand{\N}{\mathbb{N}}
\newcommand{\PP}{\mathcal{P}}
\newcommand{\F}{\mathcal{F}}
\newcommand{\Q}{\mathbb{Q}}
\newcommand{\Z}{\mathbb{Z}}

\renewcommand{\star}{{}^\ast}
\renewcommand{\phi}{\varphi}
\renewcommand{\epsilon}{\varepsilon}

\theoremstyle{definition}
\newtheorem{theorem}{Theorem}[section]

\newtheorem{remark}[theorem]{Remark}
\newtheorem{proposition}[theorem]{Proposition}

\newtheorem{definition}[theorem]{Definition}

\Crefname{question}{Question}{Questions}
\Crefname{fact}{Fact}{Facts}

\newtheoremstyle{named}{}{}{}{}{\bfseries}{.}{.5em}{\thmnote{#3}}
\theoremstyle{named}

\begin{document}

\title{A van der Waerden-free proof of Rado's theorem}
\maketitle
\author{Mauro Di Nasso  \thanks{Dipartimento di Matematica, Universit\`a di Pisa, Largo Bruno Pontecorvo 5, 56127 Pisa, Italy, email: mauro.di.nasso@unipi.it}}

\author{Lorenzo Luperi Baglini  \thanks{Dipartimento di Matematica, Università  di Milano, Via Saldini 50, 20133 Milano, Italy, email: lorenzo.luperi@unimi.it}}

\begin{abstract} 
We present a proof of the sufficiency of Rado's condition 
for the partition regularity of linear Diophantine equations that avoids 
any use of van der Waerden's theorem.
The proof is based on fundamental properties that are common knowledge 
in combinatorics of numbers and is entirely ``elementary''.
\end{abstract}

\section*{Introduction}

A fundamental result in arithmetic Ramsey Theory is \emph{Rado's Theorem}
that characterizes the linear Diophantine equations $\sum_{i=1}^nc_i x_i=0$
that are partition regular as those where
a partial sum of the coefficients equals zero. 
This latter property is known as ``Rado's condition.''
The original proof of Rado's Theorem is grounded on the observation that 
if Rado's condition holds, then configurations made of long enough arithmetic progressions and appropriate multiples of the common difference contain solutions 
to $\sum_{i=1}^{n}a_{i}x_{i}=0$; to complete the proof,
one then applies a strengthening of van der Waerden's Theorem about the existence of arbitrarily long 
monochromatic arithmetic solutions. 
To the best of our knowledge, no elementary proof of the sufficiency of 
Rado's condition avoiding any use of strengthenings of van der Waerden's Theorem has ever been given, 
where by ``elementary proof'' we mean a proof that does 
not appeal to any advanced machinery coming from other mathematical areas 
and that are often used in arithmetic number theory, 
such as ergodic theory, discrete Fourier analysis, algebra on the space
of ultrafilters, or nonstandard analysis.
This paper is aimed at filling this void. 

\smallskip
We aim to keep this paper self-contained, 
so in the first two sections we recall the few notions and properties that we need
to give our proof of Rado's Theorem, presented in Section 3. 
In the following Section 4 we develop the ``elementary approach'' and 
give an alternative proof of a result about the possibility
of combining partition regularity of different equations,
that have been previously proved by using algebra in the space of ultrafilters;
as an example of application we give a second elementary proof of Rado's Theorem.
In the final Section 5 we outline possible generalizations of our work,
and pose two open problems.

\medskip

\section{Terminology}

Let us now recall the basic notions and fix terminology.

\begin{definition}
{\rm A finite \emph{coloring} of a set $X$ is a finite partition $X=C_1\cup\ldots\cup C_r$.
The pieces $C_i$ of the partition are called \emph{colors}.

Given a finite coloring as above, a set $F$ is called \emph{monochromatic} 
if all elements of $F$ have the same color, \emph{i.e.} $F\subseteq C_i$ for some $i$.}
\end{definition}

\begin{definition}
A family $\F\subseteq\PP(X)$ is called \emph{partition regular} 
on $X$ if for every finite coloring of $X$ there exists an element $F\in\F$ which is monochromatic.

A family $\F$ is called \emph{strongly partition regular} 
if for every $A\in\F$ and for every coloring
$A=C_1\cup\ldots\cup C_r$, one of the colors $C_i\in\F$.
\end{definition}

In what follows, we will be interested in the case where $X=\N$ is the set of natural numbers; we agree that $\N$ denotes the set of positive integers,
and $\N_0=\N\cup\{0\}$ denotes the set of non-negative integers.

\smallskip
Historically, the first example of a partition regular property involving equations was provided by
a result of Schur \cite{sh} that appeared in 1916. 

\begin{theorem}[Schur]\label{schur}
For every finite coloring of $\N=C_1\cup\ldots\cup C_r$ there exists a monochromatic ``Schur triple'' 
$a, b, a+b \in C_i$.
\end{theorem}

A fundamental result in this area was then proved by van der Waerden \cite{vdw} in 1927.
With an ingenious and rather intricate double induction procedure, he proved the following.\footnote
{~A detailed proof can be found \emph{e.g.} in \cite[\S 2.1]{grs}.}

\begin{theorem}[van der Waerden]
For every $\ell\in\N$ and for every finite coloring of \, $\N=C_1\cup\ldots\cup C_r$ 
there exists a monochromatic $\ell$-term arithmetic progression
$a, a+d, \ldots, a+(\ell-1)d\in C_i$. Consequently, one color contains arbitrarily long
arithmetic progressions.
\end{theorem}

A strengthening of van der Waerden's result, that includes Schur's Theorem as a particular case, allows also a fixed multiple of the common distance of the arithmetic progression to be taken of the same color.\footnote
{~A detailed proof can be found \emph{e.g.} in \cite[\S 2.2]{grs}.}

\begin{theorem}\label{vdWB}
For every $\ell,s\in\N$ and for every coloring of \, $\N=C_1\cup\ldots\cup C_r$ 
there exists a monochromatic $\ell$-term arithmetic progression
$a, a+d, \ldots, a+(\ell-1)d\in C_i$ so that also $sd\in C_i$.
\end{theorem}

We remark that also the multiplicative versions of the above results hold,
where additive Schur's triples are replaced by multiplicative triples ``$a, b, ab$",
where arithmetic progressions are replaced by geometric
progressions $a, ad, ad^2,\ldots, ad^{\ell-1}$, and where multiples of the common difference are replaced by powers of the common ratio, respectively.

\smallskip
Since its publication, van der Waerden's result
has become one of the cornerstones of arithmetic Ramsey theory,
and a fundamental tool that is found in many proofs. 
A typical example is given by Rado's Theorem (see below)
about partition regular equations.


\begin{definition}
{\rm Let $f:\N^{k}\rightarrow\N_{0}$.
We say that the equation ``$f(x_1,\ldots,x_k)=0$''
is \emph{partition regular} on $\N$ to mean that the family of solutions
$\F_f=\{\{n_1,\ldots,n_k\}\in[\N]^k\mid f(n_1,\ldots,n_k)=0\}$ is partition regular on $\N$.
(This means that for every coloring $\N=C_1\cup\ldots\cup C_r$ there exist
monochromatic $n_1,\ldots,n_k\in C_i$ which are a solution $f(n_1,\ldots,n_k)=0$.)

In this case, by abusing terminology, sometimes we will directly say
the function $f$ is \emph{partition regular} on $\N$.}
\end{definition}

Partition regularity can be seen as a property that 
indicates the existence of a large and ``well-spread''
set of solutions, in the sense that in every given finite partition of the natural numbers
one always find elements of the same color that are a solution.

\smallskip
Note that Schur's Theorem states the partition regularity of the equation $x_1+x_2-x_3=0$.
A few years later, in 1933 R. Rado \cite{ra}  found a simple full characterization of 
all linear Diophantine equations that are partition regular.\footnote
{~A detailed historical account of the early years of Ramsey theory,
and in particular the discovery of the theorems of Schur, van der Waerden and Rado, can be found in \cite{so}.}

\begin{theorem}[Rado]
Let $a_1,\ldots,a_n\in\Z\setminus\{0\}$. Then the following are equivalent:
\begin{enumerate}
\item 
The equation $\sum_{i=1}^{n}a_{i}x_{i}=0$ is partition regular on $\N$.
\item
The following ``Rado's condition'' holds:
\\
There exists a nonempty $I\subseteq\{1,\dots,n\}$ such that $\sum_{i\in I} a_i=0$.
\end{enumerate}
\end{theorem}

The necessity of Rado's condition can be proven directly via elementary 
arithmetic arguments; see for instance \cite[\S 3.2]{grs} for details. 
As for the converse implication,
the original proof is based on the observation that a given linear equation $\sum_{i=1}^{n}a_{i}x_{i}=0$ that satisfies Rado's condition can be solved, for appropriate $s,\ell$, inside sets of the form $\{a,sd,a+d,\dots,a+\ell d\}$; the partition regularity property then
follows by an application of Theorem \ref{vdWB}. Our goal is to provide an elementary proof of the sufficiency of Rado's condition that avoids the use of van der Waerden's Theorem or its strengthenings.

\section{Preliminaries}\label{prelim}

In this section we will recall all the needed basic notions and properties
involved in our proof of the sufficiency of Rado's condition.
In order to make this paper self-contained, and to show that all the arguments we need 
are indeed elementary, we will also include most of the proofs of these properties. 

For $A\subseteq\N$ and $n\in\N_0$, 
the \emph{leftward additive shift} of $A$ by $n$
and the \emph{rightward additive shift} by $n$ are the following sets, respectively:
$$A-n=\{m\in\N\mid n+m\in A\}\,;\quad
n+A=\{n+m\in\N\mid m\in A\}.$$

The \emph{leftward} and \emph{rightward multiplicative shift} are defined similarly:
$$A/n=\{m\in\N\mid n\cdot m\in A\}\,;\quad
n\cdot A=\{n\cdot m\mid m\in A\}.$$

Shifts can be used to introduce several
elemental notions in combinatorics of numbers, including that of Banach density, sometimes called ``uniform density''.

\begin{definition}
{\rm Let $A\subseteq\N$. The \emph{upper Banach density} (or simply the \emph{Banach density}) of $A$ 
is defined by setting:
$$\text{BD}(A)=\limsup_{n\to\infty}\left(\max_{x\in\N}\frac{|(A-x)\cap[1,n]|}{n}\right)=
\limsup_{n\to\infty}\left(\max_{x\in\N}\frac{|A\cap[x+1,x+n]|}{n}\right).$$}
\end{definition}

The above limit superior can be proved to be an actual limit. 
Indeed, observe that the sequence $a_n:=\max_{x\in\N}|(A-x)\cap[1,n]|$ 
is subadditive, \emph{i.e.}, $a_{n+m}\le a_n+a_m$ for all $n<m$; 
then by Fekete's Lemma one has that $(a_n\mid n\in\N)$ converges and
$\lim_{n\to\infty}\frac{a_n}{n}=\inf_{n}\frac{a_n}{n}$.
So we have:
$$\text{BD}(A)=\lim_{n\to\infty}\left(\max_{x\in\N}\frac{|(A-x)\cap[1,n]|}{n}\right)=
\inf_{n\in\N}\left(\max_{x\in\N}\frac{|(A-x)\cap[1,n]|}{n}\right).$$

\begin{remark}
{\rm A more commonly used notion of density is that of \emph{upper asymptotic density},
which is defined by setting:
$$\overline{d}(A)=\limsup_{n\to\infty}\frac{|A\cap[1,n]|}{n}.$$
Note that $\text{BD}(A)\ge\overline{d}(A)$ since in Banach density one
considers arbitrary intervals instead of just initial segments $[1,n]$.
It is not difficult to find examples of sets $A$ such that $\overline{d}(A)=0$ but $\text{BD}(A)=1$.}
\end{remark}

Clearly $0\le\text{BD}(A)\le 1$. Note that Banach density is invariant under additive shifts;
besides, multiplicative shifts preserve non-zero density.

\begin{proposition}
Let $A\subseteq\N$ and $n\in\N$. Then

\begin{enumerate}
\item
$\text{BD}(A-n)=\text{BD}(n+A)=\text{BD}(A)$.
\item
$\text{BD}(n\cdot A)=\frac{\text{BD}(A)}{n}$.
\end{enumerate}
\end{proposition}

It is readily seen that
$\text{BD}(A\cup B)\le\text{BD}(A)+\text{BD}(B)$ for all $A,B\subseteq\N$.

Banach density is not additive; indeed one easily finds
examples of disjoint sets $A\cap B=\emptyset$ such that
$\text{BD}(A\cup B)<\text{BD}(A)+\text{BD}(B)$.
However, it behaves like a finitely additive measure in the case
of families of pairwise disjoint shifts of a given set.

\begin{proposition}\label{BDadditivity}
Let $A\subseteq\N$ and let $x_1<\ldots<x_k$.
If $(A-x_i)\cap(A-x_j)=\emptyset$ for all $i\ne j$
then $\text{BD}\left(\bigcup_{i=1}^k(A-x_i)\right)=
k\cdot\text{BD}(A)$.
\end{proposition}

\begin{proof}
Let us prove the property when $k=2$; the general case is similar.
For every $n$, pick an interval $I_n$ of length $n$ such that
$a_n:=|A\cap I_n|=\max_{x\in\N}|(A-x)\cap[1,n]|$, so that $\text{BD}(A)=\lim_{n\to\infty}\frac{a_n}{n}$.
Let $\varepsilon_n,\eta_n$ be such that $a_n=|(A-x_1)\cap I_n|+\varepsilon_n=|(A-x_2)\cap I_n|+\eta_n$.
Clearly $|\varepsilon_n|\le x_1$ and $|\eta_n|\le x_2$, so $\lim_{n\to\infty}\frac{\varepsilon_n+\eta_n}{n}=0$.
We have the following chain of inequalities:
\begin{multline*}
2\cdot\text{BD}(A)=\lim_{n\to\infty}\frac{a_n+a_n}{n}=
\lim_{n\to\infty}\left(\frac{|(A-x_1)\cap I_n|}{n}+\frac{|(A-x_2)\cap I_n|}{n}+\frac{\varepsilon_n+\eta_n}{n}\right)=
\\
=\lim_{n\to\infty}\frac{|((A-x_1)\cup(A-x_2))\cap I_n|}{n}\le
\text{BD}((A-x_1)\cup(A-x_2))\le
\text{BD}(A-x_1)+\text{BD}(A-x_2)=2\cdot\text{BD}(A).\qedhere
\end{multline*}
\end{proof}

A notion of largeness that is widely used in combinatorics of numbers is 
obtained by considering sets of differences.

\begin{definition}
{\rm The \emph{Delta-set} of a set $X\subseteq\N$ is defined as:
$$\Delta(X)=:\{x'-x\mid x<x'\ \text{in}\ X\}.$$
A set $A\subseteq\N$ is \emph{Delta-large} if $A\supseteq\Delta(X)$
for an infinite $X$.}
\end{definition}

\begin{proposition}
The family of Delta-large sets is strongly partition regular.
\end{proposition}

\begin{proof}
This is a straight application of Ramsey's Theorem.
Indeed let $A=C_1\cup\ldots\cup C_r$ be a finite coloring of
a Delta-large set $A$. Pick an infinite $X=\{x_1<x_2<\ldots\}$
such that $\Delta(X)\subseteq A$, and define the coloring
of the pairs $[\N]^2=D_1\cup\ldots\cup D_r$ by letting $\{n<m\}\in D_i\Leftrightarrow x_m-x_n\in C_i$.
By Ramsey's Theorem there exists an infinite $H\subseteq\N$ such that $[H]^2\subseteq D_i$ is monochromatic.
This means that if we let $X':=\{x_n\mid n\in H\}$ then $\Delta(X')\subseteq C_i$.
\end{proof}

\begin{definition}
{\rm A set $A\subseteq\N$ is called \emph{thick} 
if it contains arbitrarily long intervals.}
\end{definition}

We conclude this subsection by recalling two well-known facts that 
connect sets of positive Banach density, thick sets, and Delta-large sets.

\begin{proposition}\label{thick}
Let $A\subseteq\N$. 
\begin{enumerate}
\item
$\text{BD}(A)=1$ if and only if $A$ is thick.
\item
If $A$ is thick then $A$ is Delta-large.
\end{enumerate}
\end{proposition}

\begin{proof} (1). The fact that thick sets have Banach density 1 is straightforward from the definitions.
Conversely, if a set $A$ is not thick, 
pick $n\in\N$ such that $A$ does not contain $n$ consecutive natural numbers.
It is readily seen that this forces $\text{BD}(A)\leq 1-\frac{1}{n}$.

(2). We recursively define an increasing sequence $(x_n\mid n\in\N)$
such that $x_m-x_n\in A$ for all $n<m$.
Pick any $x_1\in A$. At the inductive step,
pick an interval $I=[y, y']\subseteq T$ of length $y'-y>x_n$, and let $x_{n+1}=y'$.
Then one has $x_{n+1}-x_i\in A$ for every $i\le n$ because
$y<y'-x_n\le x_{n+1}-x_i<y'$.\end{proof}

\begin{proposition}\label{Delta-set-intersection}
Let $A\subseteq\N$ have positive Banach density $\text{BD}(A)=\alpha>0$. Then $\Delta(A)\cap\Delta(X)\ne\emptyset$
for all infinite $X\subseteq\N$.
\end{proposition}

\begin{proof}
Consider elements $x_1<x_2<\ldots<x_n$ in $X$ where $n>1/\alpha$.
If we show that there exist $i<j\le n$ with $(A-x_i)\cap(A-x_j)\ne\emptyset$ then
we are done. Indeed in this case we can pick $y\in(A-x_i)\cap(A-x_j)$, and
if $a,a'\in A$ are such that $y=a-x_i=a'-x_j$ then $x_j-x_i=a'-a\in\Delta(A)\cap\Delta(X)$.
So, assume towards a contradiction that $(A-x_i)\cap(A-x_j)=\emptyset$ for all $i<j\le n$.
Then, by Proposition \ref{BDadditivity} we would have
$\text{BD}(\bigcup_{i=1}^n(A-x_i))=n\cdot\alpha>1$.
\end{proof}

\section{Proof of Rado's Theorem}\label{proof}

We are to give an elementary proof of Rado's Theorem based
on the facts that we recalled above. 

For equations in two variables, Rado's condition reduces to the fact that 
the polynomial $P(x_1,x_2)$ is of the form
$c(x_1-x_2)$ for some $c\in\Z$.\footnote
{~This characterization can be extended to the partition regularity of arbitrary 
nonlinear equations $P(x_{1},x_{2})=0$: in fact such equations are partition 
regular if and only if they are multiples of $x_{1}-x_{2}$, as was proved in \cite{lba}.} 
In this case, the partition regularity is trivial by considering any constant solution
$P(a,a)=0$. 
 
Similarly, also polynomials $P(x_1,\ldots,x_n)$ where
the sum of all coefficients $\sum_{i=1}^nc_i=0$ admit constant solutions $P(a,\ldots,a)=0$,
and hence they are trivially partition regular.

The following proposition deals with the first non-trivial case in which
there are three variables and where the sum of all coefficients may be different from 0.

\begin{proposition}\label{3variables}
For every $c\in\N$ and for every $d\in\Z$, the Diophantine equation $cx-cy-dz=0$ 
is partition regular on $\N$. 
\end{proposition}

\begin{proof} As already observed above, if $d=0$ the result is trivial. 
If $d\neq 0$ then, by reordering the terms of the equation and by 
exchanging the variables $x$ and $y$ if necessary,
we can assume without loss of generality that $d\in\N$. 
Let $\N=C_1\cup\ldots\cup C_s\cup C_{s+1}\cup\ldots\cup C_r$ be a finite coloring of 
$\N$, where we arranged colors so that 
the Banach density $\text{BD}(C_i)>0$ if and only if $i\le s$.
The set $D=:C_1\cup\ldots\cup C_s$ has Banach density $1$, 
and so it is Delta-large by Proposition \ref{thick}.
Since Delta-large sets are strongly partition regular, there exists $i\le s$ and
an infinite $W$ such that $\Delta(W)\subseteq C_i$.
Now, $\text{BD}(C_i)>0$ implies that $\text{BD}(c\cdot C_i)>0$, and so
$\Delta(c\cdot C_i)\cap\Delta(d\cdot W)\ne\emptyset$, by Proposition \ref{Delta-set-intersection}.
This means that there exist $a_1>a_2$ in $C_i$ and $a'_1>a'_2$ in $W$ such that 
$ca_1-ca_2=da'_1-da'_2=da_3$ where $a_3=:a'_1-a'_2\in\Delta(W)$. 
Finally, observe that $a_1,a_2,a_3\in C_i$ are monochromatic. 
\end{proof}

\begin{remark}
{\rm The idea of using density arguments combined with the strong partition regularity of
Delta-sets is certainly not new. To the authors' knowledge, the first use of such arguments 
goes back to Bergelson's paper \cite{be} of 1986, where he proved
a density version of Schur's Theorem.}
\end{remark}

We can now fulfill our original goal.

\begin{theorem}\label{finale}
If a linear homogeneous polynomial $P\in\Z[x_1,\ldots,x_k]$ 
satisfies Rado's condition, then it is partition regular on $\N$.
\end{theorem}

\begin{proof}
By renaming the variables if necessary, without loss of generality we can assume 
that $P(x_1,\ldots,x_k)=\sum_{i=1}^{k}c_i x_i$ 
where $\sum_{i=1}^m c_i=0$ for some $1\le m\le k$. If $m=k$ the equation admits constant solutions, so we restrict to the case $m<k$.

Let $I=\{1,\dots,m\}$, let $I_{1}=\{i\in \{1,\dots,m\}\mid c_{i}>0\}$, and let $I_{2}=I\setminus I_{1}$. 
Observe that $\sum_{i\in I_{1}} c_{i}=-\sum_{i\in I_{2}} c_{i}$ by Rado's condition. 
Let $c=\sum_{i\in I_{1}} c_{i}, d=-\sum_{i=m+1}^{k} c_{i}$. 
Given any finite coloring of $\N$, by Proposition \ref{3variables} 
we can find monochromatic $y_{1},y_{2},y_{3}$ so that $cy_{1}-cy_{2}-dy_{3}=0$. 
By letting $x_{i}=y_{1}$ for $i\in I_{1}, x_{i}=y_{2}$ for $i\in I_{2}$ and $x_{i}=y_{3}$ 
for $i\in \{m+1,\dots,k\}$, we obtain a monochromatic solution of $P\left(x_{1},\dots,x_{k}\right)=0$.
\end{proof}

\begin{remark}
{\rm
In the previous proof, the monochromatic solutions we obtain are not constant, 
but neither are they pairwise different.
In fact, they have at least two different values, 
since we can choose $y_1>y_2$, but it could happen that $y_3=y_1$ or $y_3=y_2$.
We observe that this is already enough to obtain a proof of the
3-term van der Waerden's Theorem, because 
we have monochromatic solutions to the equation $x+y-2z=0$ where $x>y$,
and hence $y<z<x$ form a non-trivial arithmetic progression.}
\end{remark}

\section{Piecewise syndetic sets and partition regularity}

In this section we will show an example that gives further evidence
to the fact that the largeness notion of piecewise syndeticity 
for sets of integers is instrumental in partition regularity problems.
Precisely, we will consider a strong property for combining partition
regularity of different equations. 
While its known proof \cite{lb}
heavily relies on the ultrafilter space $\beta\N$ and its algebra,
here we will show that an alternative proof can be given
that is only based on the fundamental notions  
introduced so far, and is entirely ``elementary''.
At the end of the section,
we will also show an alternative proof of Rado's Theorem
based on the above mentioned property.

\subsection{Thickness and syndeticity}
We already introduced the notion of (additively) thick set,
as a set that contains arbitrarily long intervals.
The multiplicative counterpart is the following.

\begin{definition}
A set $A\subseteq\N$ is \emph{multiplicatively thick} if it contains arbitrarily large
sets of multiples.
\end{definition}


Recall the following elemental notions of combinatorics of numbers.

\begin{definition}
\

\begin{itemize}
\item 
$A$ is \emph{additively syndetic} if it has ``bounded gaps'',
\emph{i.e.} if there exists $k\in\N$ such that for every interval $I$
of length $|I|\ge k$ one has that $A\cap I\ne\emptyset$.

\item 
$A$ is \emph{multiplicatively syndetic} if it has ``bounded gaps on multiples'',
\emph{i.e.} if there exists $k\in\N$ such that for every $x\in\N$ 
one has that $A\cap\{x, x\cdot 2,\ldots,x\cdot k\}\ne\emptyset$.

\item 
$A$ is \emph{additively piecewise syndetic} 
if it is the intersection of an additively thick set with
an additive syndetic set.

\item 
$A$ is \emph{multiplicatively piecewise syndetic} 
if it is the intersection of a multiplicatively thick set with
a multiplicatively syndetic set.
\end{itemize}
\end{definition}

We itemize below the first properties, that are basic knowledge in this area.
Proofs are directly obtained from the definitions and are omitted.

\begin{proposition}
Let $A\subseteq\N$. Then

\begin{enumerate}
\item
$A$ is additively syndetic if and only the complement $A^c$ is not additively thick
if and only if $\N=\bigcup_{i=1}^\ell (A-n_i)$ is a finite union of
leftward additive shifts of $A$.

\item
$A$ is multiplicatively syndetic if and only the complement $A^c$ is not multiplicatively thick
if and only if $\N=\bigcup_{i=1}^\ell A/n_i$ is a finite union of
leftward multiplicative shifts of $A$.
\end{enumerate}
\end{proposition}

A well-known relevant property of piecewise syndetic sets that is neither satisfied
by thick nor by syndetic sets, is partition regularity.

\begin{proposition}\label{PRofPS} 
The family of additively (or multiplicatively) piecewise syndetic sets
is strongly partition regular.
\end{proposition}

\begin{proof} 
Below, we do not specify whether we are considering the additive or the multiplicative notion, 
since the arguments used in the proof are exactly the same for both cases.

Let $A\subseteq\N$ be piecewise syndetic.
By induction on $r$ we will show that in every
partition $A=C_1\cup\ldots\cup C_r$ one of the pieces $C_i$ is piecewise syndetic.
The base case $r=1$ is trivial.
At the inductive step $r>1$, let $C':=C_1\cup\ldots\cup C_{r-1}$.
We note that, using the inductive hypothesis, it suffices to show that one 
of $C'$ and $C_r$ is piecewise syndetic.

Pick a syndetic set $S$ and a thick set $T$ such that $A=S\cap T=C'\cup C_r$.
We distinguish two cases.
If $S\setminus C'$ is syndetic, then 
$C_r=(S\setminus C')\cap T$ is piecewise syndetic as the intersection of a syndetic and a thick set;
if instead $S\setminus C'$ is not syndetic then the complement 
$(S\setminus C')^c=S^c\cup C'$ is thick, and in this case
$C'=S\cap(S^c\cup C')$ is piecewise syndetic.
\end{proof}

\begin{remark}
{\rm The previous result is a particular case of a more general property.
Precisely, along the same lines, one can prove the following:
\begin{itemize}
\item
\emph{Let $\F\subseteq\mathcal{P}(X)$ be a nonempty family 
that is upward closed with respect to inclusion. Then the 
following family is partition regular on $X$}:\footnote
{~The family $\F^*=\{B^c\mid B\in\F\}$ is usually called the \emph{dual family} of $\F$.}
$$\mathcal{F}':=\{A\cap B\mid A\in\mathcal{F},\ B^c\notin\mathcal{F}\}$$
\end{itemize}}
\end{remark}

\subsection{The principle of compactness}

\smallskip
A well-known tool in combinatorics is given by the compactness principle,
a general property that connects infinitary and finitary properties.
There are several different formulations of that fundamental principle;
the one we will use in the sequel is the following.

\begin{theorem}[Compactness principle]\label{compactness}
Suppose that $\F\subseteq\PP(X)$ only contains finite sets.
Then for every $r\in\N$,
the family $\F$ is $r$-partition regular on $X$ if and only 
there exists a finite $Y\subseteq X$
such that $\F\cap\PP(Y)$ is $r$-partition regular on $Y$.
\end{theorem}

By ``$r$-partition regularity'' we mean that the partition regularity property
holds for $r$-colorings, \emph{i.e.}, for colorings into $r$-many colors.

\begin{proof} One implication is immediate. To prove that the other implication is ``elementary'', we include its proof here. The argument is the same as, \emph{e.g.},
in \cite[\S 1.5]{grs}. 

Consider the topological space of $r$-colorings $\{1,\dots,r\}^X$ 
obtained as the topological product of the discrete space $\{1,\ldots,r\}$.
By Tychonoff's Theorem, such a space is compact.\footnote
{~Note that here we are using the \emph{axiom of choice}.}
Now assume for the sake of contradiction that
that for every finite $Y\subseteq X$, the family $\F\cap\PP(Y)$ is 
not $r$-partition regular on $Y$, \emph{i.e.},
the set $C_Y$ of the $r$-colorings 
such that no monochromatic $F\in\F$ is contained in $Y$ is nonempty. 
Now observe that $C_{Y_1\cup\ldots\cup Y_n}\subseteq C_{Y_1}\cap\ldots\cap C_{Y_n}$,
and so the family $\{C_Y\mid Y\subseteq X\ \text{finite}\}$ has the finite intersection property.
Since all sets $C_Y$ are closed, by compactness
the intersection $C_{X}=\bigcap_{Y\subseteq X, Y \ \text{finite}} C_{Y}$
is nonempty, contradicting the $r$-partition regularity of $\F$ on $X$.
\end{proof}

For example, by compactness one has that 
Schur's Theorem \ref{schur} implies the following finitary statement 
(the converse implication is trivial):\footnote
{~Actually, both the original Schur's Theorem \cite{sh}
and van der Waerden's Theorem were stated in their finitary versions.}
\begin{itemize}
\item
\emph{For every $r\in\N$ there exists $n\in\N$ such that for every
$r$-coloring $\{1,\ldots,n\}=C_1\cup\ldots\cup C_r$ there exists
a monochromatic triple $a,b,a+b\in C_i$.}
\end{itemize}

\subsection{Combining partition regularity of different equations}

Let us start with a disclaimer: all the results we will prove in this subsection 
(although sometimes formulated slightly differently) have already been 
proved in \cite{lb} and \cite{dlb}.\footnote
{~Precisely, Proposition \ref{shiftPS} is a straight consequence of Theorem 4.13 of \cite{lb}, 
Proposition \ref{PR-MPS} was proven in Corollary 6.14 of \cite{lb}, 
Proposition \ref{jointPR} was proven in Lemma 2.1 of \cite{dlb}, 
Lemma \ref{inductive-lemma} (in a much-strengthened way) 
can be found in Theorem 2.7 of \cite{dlb}, and the argument in the proof of
Theorem \ref{finale2} is similar to the one used to prove 
Corollary 2.5 and Theorem 2.11 in \cite{dlb}.}
However, all the proofs there were carried within the ultrafilter space $\beta\N$
and used in an essential way its algebraic properties; instead,
here we will provide alternative proofs that use only the elementary tools
introduced so far.

The crucial tool we will need is a sort of universal property of the piecewise syndetic
sets with respect to partition regularity.
To the author's knowledge, it was first proven in \cite{mc} (see Theorem 1.2)
by using compactness of the topological spaces $\{1,\ldots,r\}^\N$, 
but has since been considered several times in the literature; 
for example it can be obtained as a direct consequence of Theorem 4.13 of \cite{lb}, 
which was proven using the ultrafilter method.
We give here a short proof as a straight application of combinatorial compactness.

\begin{proposition}\label{shiftPS}
Let $\F\subseteq\PP(\N) $ be a partition regular family of finite sets.
Then every additive (or multiplicative) piecewise syndetic set contains an 
additive (or multiplicative, resp.) shift of an element of $\F$.
\end{proposition}

\begin{proof} 
We consider here only the multiplicative case, 
because the additive one is proved by using exactly the same arguments.

Given a multiplicatively piecewise syndetic set $A$,
let $r\in\N$ be such that the union $T:=\bigcup_{i=1}^r A/i$ is multiplicatively thick. 
Since $\F$ only contains finite sets, by compactness 
there exists $n$ such that for every $r$-coloring of $[1,n]$ 
there exists a monochromatic set $F\in\F\cap\PP([1,n])$. 
By multiplicative thickness, 
there exists $m$ such that the multiplicative shift $m\cdot [1,n]\subseteq T$.
Now consider the $r$-coloring $[1,n]=C_1\cup\ldots\cup C_r$
where for every $x\in [1,n]$ one lets $x\in C_j$ if and only if $j=\min\{i\leq k\mid m x\in A/i\}$. 
Let $F\in\F\cap \PP([1,n])$ be monochromatic, say $F\subseteq C_i$.
Then, by definition, for every $x\in F$ we have $mx\in A/i$, 
and hence the multiplicative shift $(im)\cdot F\subseteq A$, as desired.
\end{proof}

\begin{remark}
{\rm With essentially the same proof, the result above can be extended
to any semigroup $(S,\star)$ in place of $(\N,+)$ or $(\N,\cdot)$.
Some caution is needed in the non-commutative case,
where one needs to consider the one-sided notions of 
right (or left) thick set, and of right (or left, resp.) piecewise syndetic set.\footnote
{~See \S 2 of \cite{bhm} or \S 4.4 of \cite{hs}.}}
\end{remark}

The central role of multiplicatively piecewise syndetic sets
for our purposes is shown by the following property.
We prove it as a direct consequence of the 
universal property of piecewise syndetic sets that we have seen above.

\begin{proposition}\label{PR-MPS}
Let $P\in\Z[x_1,\ldots,x_k]$ be a homogeneous polynomial. 
Then the following are equivalent:
\begin{enumerate}
\item 
The equation $P(x_1,\ldots,x_k)=0$ is partition regular on $\N$.
\item 
For every multiplicatively piecewise syndetic set $A\subseteq\N$ 
there exist elements $a_1,\ldots,a_k\in A$ that are a solution, \emph{i.e.},
$P(a_1,\ldots,a_k)=0$.
\end{enumerate}
\end{proposition}

\begin{proof} 
$(2)\Rightarrow(1)$ directly follows from the partition regularity
of multiplicatively piecewise syndetic sets (see Proposition \ref{PRofPS}).

\smallskip
$(1)\Rightarrow (2)$. Consider the family
$$\F=:\left\{\{a_1,\ldots,a_k\}\in[\N]^k\mid P(a_1,\ldots,a_k)=0\right\}.$$
By the hypothesis, $\F$ is partition regular on $\N$ and so, by the previous Proposition \ref{shiftPS},
there exists a set $\{b_1,\ldots,b_k\}\in\F$ and $y\in\N$ 
such that the multiplicative shift
$y\cdot\{b_1,\ldots, b_k\}\subseteq A$.
If we let $a_i:=y\cdot b_i\in A$ then we have
$P(a_1,\ldots,a_k)=y^\ell\cdot P(b_1,\ldots,b_k)=0$,
where $\ell$ is the degree of homogeneity of $P$.
\end{proof}

Combining the two previous properties, one
obtains a valuable tool that allows combining different 
equations in such a way as to preserve partition regularity.

\begin{theorem}\label{jointPR}
Suppose the homogeneous polynomials $P_i\left(x_{i,1},\ldots,x_{i,k_i}\right)$ 
are partition regular on $\N$ for $i=1,\ldots, N$. Then in every multiplicatively piecewise syndetic set $A$
one finds elements $b$ and $a_{i,2},\ldots,a_{i,k_i}$
such that $P_i\left(b,a_{i,2},\ldots,a_{i,k_i}\right)=0$ for $i=1,\ldots,N$.
\end{theorem}

\begin{proof}
For $i=1,\ldots,N$, consider the set:
$$\Gamma_i:=\{y\in A\mid \forall y_2,\ldots,y_{n_i}\in A\  P_i(y,y_2,\ldots,y_{n_i})\ne 0\}.$$
We observe that $\Gamma_i$ is not multiplicatively piecewise syndetic because it does not contain
any solution of $P_i$, and therefore neither is the finite union 
$\Gamma:=\bigcup_{i=1}^N\Gamma_i\subseteq A$.
By strong partition regularity, it follows that
$A\setminus\Gamma$ is multiplicatively piecewise syndetic, and hence nonempty. Finally, by definition,
if $b$ is any element of $A\setminus\Gamma$ then
for every $i=1,\ldots,N$ there exist elements $a_{i,2},\ldots,a_{i,n_i}\in A$ 
such that $P_i(b,a_{i,2},\ldots,a_{i,n_i})=0$, as desired. 
\end{proof}

As a consequence of the previous theorem we obtain the following general property
about partition regularity.

\begin{proposition}\label{inductive-lemma}
Assume the homogeneous polynomial $P\in\Z[x_1,\ldots,x_k]$ is partition regular on $\N$.
Then for every $q\in\Q\setminus\{0\}$ also the polynomial
$$P_q(x_1,\ldots,x_{k-1},x_k,x_{k+1})=P(x_1,\ldots,x_{k-1},x_k+qx_{k+1})$$
is partition regular on $\N$.
\end{proposition}

\begin{proof}
Let $A$ be a multiplicatively piecewise syndetic, and let $q=\frac{d}{c}$ 
where $c\in\N$ and $d\in\Z\setminus\{0\}$.
By Proposition \ref{3variables}, the 
homogeneous polynomial $Q(y_1,y_2,y_3):=cy_1-cy_2-dy_3$ is partition regular. 
Now apply Proposition \ref{jointPR} to the polynomials $P$ and $Q$, 
and obtain the existence of elements $b$ and $a_1,\ldots,a_{k-1},b_2,b_3$ in $A$ such that
$P(a_1,\ldots,a_{k-1},b)=0$ and $Q(b,b_2,b_3)=c b-c b_2-db_3=0$.
Since $b=b_2+q b_3$, the elements
$a_1,\ldots,a_{k-1},b_2,b_2+q b_3\in A$ are a solution of 
$P_q$, as desired. 
\end{proof}

Several applications of the above theorem and its consequences
can be found in \cite{dlb}.
Here, to illustrate its use, we give a second, alternative 
elementary proof of the sufficiency of Rado's condition. 

\begin{theorem}\label{finale2}
If the linear homogeneous polynomial $P\in\Z[x_1,\ldots,x_k]$ 
satisfies Rado's condition, then it is partition regular on $\N$.
\end{theorem}

\begin{proof}
By renaming the variables if necessary, without loss of generality we can assume 
that $P(x_1,\ldots,x_k)=\sum_{i=1}^{k}c_i x_i$ 
where $\sum_{i=1}^m c_i=0$ for some $1\le m\le k$.
The polynomial $Q(x_1,\ldots,x_m)=:c_1x_1+\ldots+c_m x_m$ is trivially partition regular, 
as it is solved by every constant $m$-tuple.

Then, by the previous proposition, also the following polynomial is partition regular:
$$Q\left(x_1,\ldots,x_{m-1},x_m+\frac{c_{m+1}}{c_m}x_{m+1}\right)=
c_1x_1+\ldots+c_m x_m+c_{m+1}x_{m+1}.$$
By iterating the procedure for $k-m$ steps, we finally
obtain that $P$ is partition regular.
\end{proof}

\medskip
\section{Conclusions}\label{conclusions}

The version of Rado's Theorem considered here is a special case of 
the general formulation that characterizes the partition regularity on $\N$ of 
arbitrary (finite) systems of Diophantine equations. 
We have not been able to adapt our methods to the general formulation 
of Rado's Theorem, hence our first question is the following.

\medskip
{\bfseries Question 1:} 
Is it possible to give an elementary proof of the general Rado's Theorem 
for systems of linear equations, without resorting to Van der Waerden's Theorem?

\smallskip

In our proof of Theorem \ref{finale}, the monochromatic solutions we produce are ``almost"
constant'', in the sense that they attain only 3 different values (at most). 
Often, in Ramsey theory we prefer to avoid such restrictions. Hence, our final question is the following.

\medskip
{\bfseries Question 2:} 
Is it possible to give an alternative elementary proof of Rado's Theorem that proves also the 
existence of injective monochromatic solutions
(\emph{i.e.}, solutions whose coordinates are pairwise distinct)?

\bigskip
\noindent
\textbf{Acknowledgements.}
The authors thank D\"om\"ot\"or P\'alv\"olgyi for suggesting a simplification
in our original proof of Theorem \ref{finale}.
The authors also thank the referee for carefully reading the paper
and for the useful indications.

Both authors  were supported by the project PRIN 2022 ``Logical methods
in combinatorics'', 2022BXH4R5, Italian Ministry of University and Research (MUR).

\section*{Conflicts of interest}

The authors do not work for, advise, own shares in, or receive funds from any organization
that could benefit from this article, and have declared no affiliations other than their research
organizations.

\end{document}